\documentclass{amsproc}
\usepackage{amsmath,amssymb,graphicx}

\usepackage{pdfsync}

\theoremstyle{definition}

\theoremstyle{remark}

\numberwithin{equation}{section}

\newcommand{\R}{{\mathbb R}}

\newcommand{\be}{\begin{equation}}
\newcommand{\ee}{\end{equation}}

\begin{document}

\title[Transverse stability in KP equations]{Transverse stability of periodic traveling waves in Kadomtsev-Petviashvili equations: A numerical study}

\author[C. Klein]{Christian Klein}
 \address{Institut de Math\'ematiques de Bourgogne \\ 9 avenue Alain Savary\\
21078 Dijon Cedex} 
 \email{christian.klein@u-bourgogne.fr}
\thanks{The first author was supported by the project FroM-PDE funded by the European
Research Council through the Advanced Investigator Grant Scheme, the Conseil R\'egional de Bourgogne
via a FABER grant and the ANR via the program ANR-09-BLAN-0117-01.}

\author{Christof Sparber}
\address{Department of Mathematics, Statistics, and Computer Science,
University of Illinois at Chicago,
851 South Morgan Street
Chicago, Illinois 60607, USA}
\email{sparber@uic.edu}
\thanks{The second author was supported by the Royal Society through is University Research Fellowship.}

%    General info
\date{\today}


\begin{abstract}
We numerically investigate transverse stability and instability of 
so-called cnoidal waves, i.e., periodic traveling wave solutions of 
the Korteweg-de Vries equation, 
under the time-evolution of the Kadomtsev-Petviashvili equation. In particular, we find that in KP-I  small amplitude cnoidal waves are stable 
(at least for spatially localized perturbations) and only become unstable above a certain threshold. In contrast to that, KP-II is found to be stable for all 
amplitudes, or, equivalently, wave speeds. 
This is in accordance with recent analytical results for solitary waves given in \cite{RT1, RT2}. 
\end{abstract}

\maketitle

\section{Introduction}

The starting point of this investigation is the celebrated {\it 
Korteweg-de~Vries equation} (KdV) for shallow water waves \cite{KdV}, i.e. 
\be
\partial_t u +   u \, \partial_x u +   \partial_{xxx} u   =0 ,\quad t, x \in \R,
\label{KdV}
\ee
subject to some given initial data $u|_{t=0}=u_0(x) \in \R$.
Clearly, \eqref{KdV} only accounts for uni-directional waves, where in fact one already has taken into account the change of coordinate $x\mapsto x+ t$ in order 
to transform the arising linear transport operator $$T_\pm 
u:=\partial_t u + \partial_x u$$ into a single partial derivative w.r.t.\ time, only. 
The corresponding two-dimensional analogue for waves with finite speed of propagation is given by $\partial_t u +  \sqrt{-\Delta}u $. Here, $ \sqrt{-\Delta}$ 
is the Fourier multiplier corresponding to the symbol $\sqrt{\xi_x^2 + \xi_y^2}$, where $\xi_x, \xi_y\in \R$ denote the frequencies (or Fourier variables) in 
$x$ and $y$ direction, respectively. A generic approach for taking into account \emph{weak transverse effects} in the propagation of quasi uni-directional waves is as follows: By Taylor expanding the dispersion relation
$$
 \sqrt{ \xi_x^2 + \xi_y^2} \approx \xi_x + \frac{1}{2} \xi^{-1}_x \xi_y^2, \quad \text{when $| \xi_x | \ll 1$ and $ \Big |\frac{\xi_y}{\xi_x} \Big | \ll 1$},
$$
and taking into account the same change of reference frame as before, one formally obtains
$$
\partial_t u + \sqrt{-\Delta}u  \approx \partial_t u + \frac{1}{2} \partial_x^{-1} \partial_{yy} u,
$$
where $\partial_x^{-1} \partial_{yy} $ is defined as a singular Fourier multiplier with symbol $-i\xi^{-1}_x \xi_y^2$. 
By following this approach one formally derives the 
well known \emph{Kadomtsev-Petviashvili equation} \cite{KP} (KP), which after rescaling $y \mapsto \sqrt{2} y$, can be written as:
\be \label{KP}
\partial_t u +   u \, \partial_x u  +  \partial_{xxx} u   +\lambda  \partial_x^{-1} \partial_{yy} u = 0 , \quad t,x,y, \in \R, \ \lambda = \pm 1, 
\ee
subject to initial data $u|_{t=0} = u_{\rm in} (x,y)$, such that $u_{\rm in}(x,0) = u_0(x)$.
Equation \eqref{KP} can be seen as a (generic) $2+1$ dimensional generalization of 
the KdV equation for quasi-uni-directional waves with weak transverse 
effects (see also \cite{GS}). The case $\lambda = +1$ corresponds to the so-called KP-II (or defocusing 
KP) equation, valid for strong surface tension, whereas $\lambda = -1$ corresponds to the KP-I (or focusing KP) equation, obtained in the case of small surface tension. 
Two basic conservation laws for \eqref{KP} are the {\it mass}
$$
M = \iint_{\R^2}  u^2 dx \,dy ,
$$
and the {\it energy}
$$
E= \iint_{\R^2} (\partial_x u)^2 - \frac{1}{3} u^3 + \lambda (\partial_x^{-1} \partial_{y} u)^2 dx\, dy.
$$

The KP equation arises, e.g., as an asymptotic model for the 
propagation of capillary-gravity waves \cite{ASL}, or for traveling waves in Bose-Einstein condensates, cf. \cite{CR}.  
Analogously to KdV, the KP equation can be written in the form of a Lax pair \cite{Z} and explicit solutions can be obtained by invoking the {\it inverse scattering transform}, see, e.g., \cite{APS, FS}. 
Generalizations to higher order nonlinearities of the form $\partial_x u^p$ are also considered (see \cite{Sa} for a broad review), with $p=4/3$ being the critical case (in the sense that the potential energy 
can be controlled by the dispersion). Thus, in general, one can expect finite-time blow-up for $p>4/3$, cf. \cite{Li} for more details. 
Clearly, any $y$-independent solution of \eqref{KP} is a solution of \eqref{KdV}, 
forming the so-called \emph{KdV sector}. It is then a natural and physically relevant question to ask whether traveling wave solutions of KdV 
are stable when considered as solution of the KP equation. This problem is usually referred to as {\it transverse (in-)stability}, to be considered in the following for so-called {\it cnoidal waves}.

\section{Cnoidal waves}

\begin{figure}[htb]
    \begin{center}
      \includegraphics[width=0.6\textwidth]{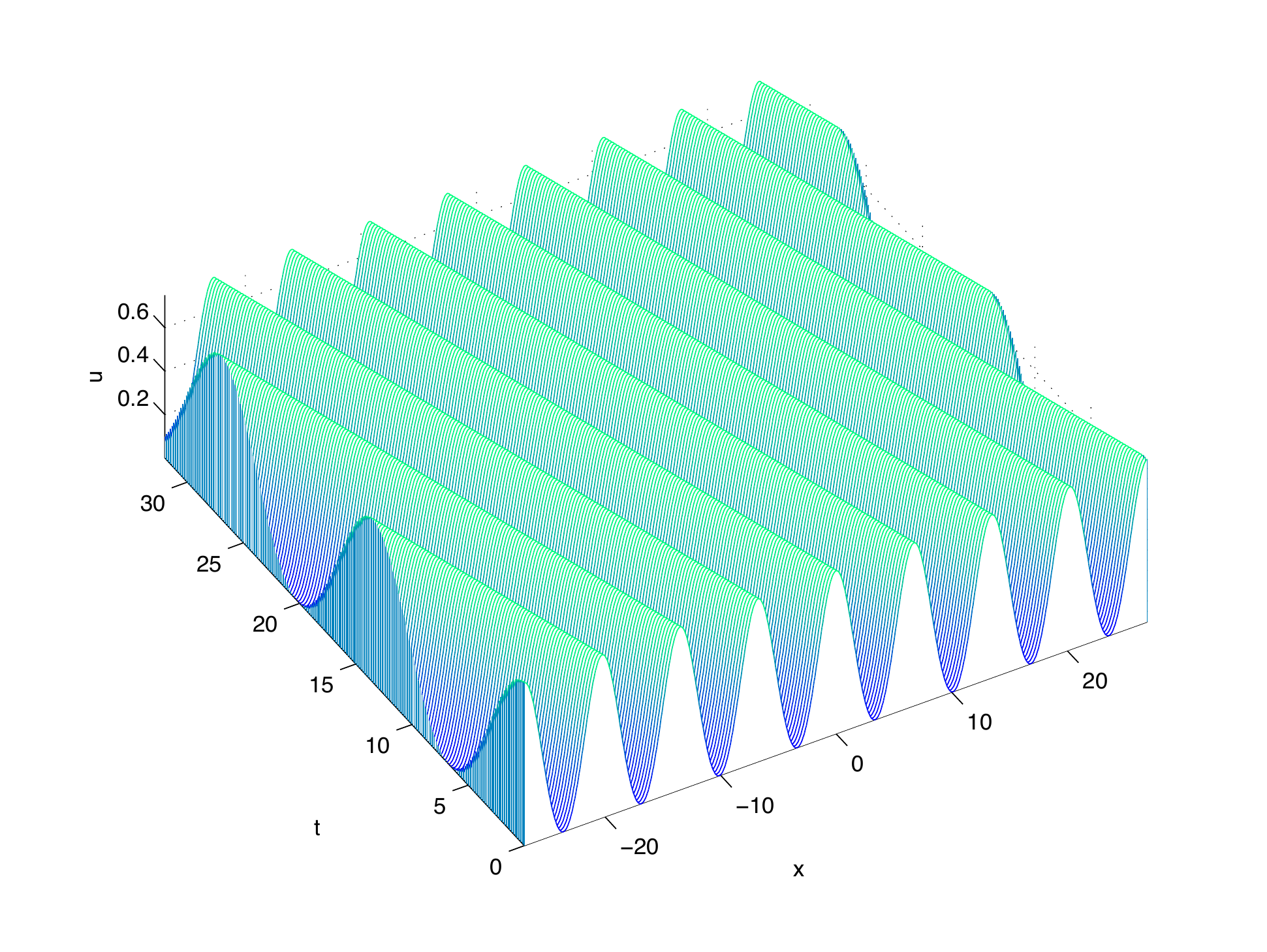}
    \end{center}
\caption{Cnoidal solution (\ref{cn}) to the KdV equation  with parameters
$\kappa=0.5$, $k=0.5$, and $x_0=u_{0}=0$.}
\label{cnoidal}
\end{figure}

Periodic traveling wave solutions to the KdV equation can be written in terms of elliptic functions. In particular the well-known cnoidal waves are given by
\begin{equation} 
    u_{\rm cn}(x,t) = u_{0}+12\kappa^{2}k^{2}\mbox{cn}^{2}(k(x-x_0 - (V+u_0)t);k), \quad V = 
    4\kappa^{2}(2k^{2}-1)
    \label{cn},
\end{equation}
where $u_{0}, \kappa , x_0 \in \R$ are arbitrary constants. Here, cn$(\cdot\, ; k)$ denotes the Jacobi elliptic cosine function \cite{BF} with elliptic modulus $k\in [0,1[$.
Each solution given by \eqref{cn} is {\it periodic} in $x$ with period 
$$
\omega(k) = \frac{2 K(k)}{\kappa} = \frac{2}{\kappa} \int_0^{\pi /2} \frac{ds}{\sqrt{1-k^2\sin^2 s}}\, , 
$$
where $K(k)$ is the complete elliptic integral of the first kind. The 
solution for $x_0=u_{0}=0$ and $\kappa=0.5$ can be seen  in 
Fig.~\ref{cnoidal}.

Note that in the limit as $k\to1$ we have: $K(k)\to +\infty$ and $\mbox{cn}(x;k)\to \mbox{sech}(x)$ in which case the cnoidal solution reduces to the well-known KdV soliton: 
\begin{equation} 
    u_{\infty }(x,t) = u_{0}+12\kappa^{2}\mbox{sech}^{2}(x-x_0 - (4 \kappa^2 + u_0) t).
\end{equation}
In contrast to cnoidal waves the KdV soliton is spatially localized, since it exponentially decays as $x\to \pm \infty$. 
Also note that the amplitude and the wave speed $V$ are correlated through the parameter $\kappa$. In the case of higher order nonlinearities, similar 
periodic traveling wave solutions appear, such as {\it dnoidal waves} (for modified KdV) or {\it snoidal waves}. All of them received their name from the fact that they can be expressed in terms of 
Jacobi elliptic functions.

There has been a large amount of work aimed at understanding the 
stability of traveling wave solutions.
The stability of solitary waves under the KdV flow is well known since the works of \cite{Be, Bo}. For periodic solutions, such as cnoidal waves, the situation is less complete and 
there are mainly two types of stability results: Orbital (or nonlinear) stability with respect to periodic and sub-harmonic perturbations has been shown in \cite{ABS, KaDe}, 
and linear and/or spectral stability with respect to non-periodic 
(i.e., localized or bounded) perturbations is proved in \cite{BoDe, HaKa}.

This situation becomes  even more involved when one  
asks about the transverse stability of traveling waves under the KP flow: Early results in \cite{APS} establish linear 
instability of the KdV soliton for KP-I, whereas for 
KP-II these solitary waves are shown to be linearly stable. More recently, rigorous results given in \cite{RT1, RT2} have established the existence of 
a {\it critical threshold} in $\kappa$ for KP-I posed on $\R_x\times \mathbb T_y$, where $\mathbb T = \R/(2\pi \mathbb Z)$. In particular, for $\kappa> \kappa_{c}=3^{-1/4}$, solitary waves are shown to 
be orbitally unstable for KP-I \cite{RT1}, whereas for $\kappa < \kappa_c$ they are proved to be orbitally (nonlinearly) stable \cite{RT2}. Generalizations of both 
results to higher order nonlinearities are possible. In addition, the full nonlinear stability (including asymptotic stability) 
of KP-II posed on $ \R_x\times \mathbb T_y$ for any $\kappa\in \R$ has been obtained in \cite{MT}. 
Note that if $u(t,x,y)$ is solution of \eqref{KP}, then so is $u(-t,-x, y)$ and hence instabilities arising for positive times, also arise for negative times. 
Numerical studies of the stability of the KdV soliton for KP have 
been performed for instance in \cite{KS}. 

Concerning the transverse (in-)stability of cnoidal waves under the KP flow, several results have been obtained in recent years: Spectral instability for small periodic solutions to 
KP-I posed on $\R_x\times \R_y$ has been proved in \cite{Ha}, where the author allows for 
two-dimensional perturbations which are either periodic in the direction of propagation (with the same period as the one-dimensional traveling wave), or non-periodic perturbations. 
In addition, spectral stability for small solutions to KP-II on $\R_x\times \R_y$ is also proved in \cite{Ha}, assuming perturbations which are  periodic in $x$ and have long wave-length in $y$.
For generalized KP-I type equations on $\R_x\times \R_y$, linear instability is proved in \cite{JZ}. Finally, \cite{HSS} proves spectral and 
linear instability for KP-I posed on $\mathbb T_x\times \mathbb T_y$, assuming that the perturbation 
admits the same period in $x$ as the cnoidal wave.

In order to establish a more qualitative picture, we shall in the following perform a numerical study of the transverse stability of cnoidal waves under the KP flow. The 
perturbations considered, are either spatially localized in both $x$ and $y$ (see Section \ref{sec:gauss}), or $y$-periodic, and thus non-localized, deformations of cnoidal waves (see Section \ref{sec:per}).  
Before describing our finding, we shall first describe the numerical algorithm in the upcoming section.

\section{Numerical Method}

In the present work, the KP equations (\ref{KP}) will be numerically studied for $(x,y,t)\in \mathbb{T}^2 \times 
\mathbb{R}$ (note that this is the same setting as used in \cite{HSS} and still close to the one of \cite{RT1, RT2}, where $y\in \mathbb T$.).
The periodicity in both spatial coordinates allows the use of {\it Fourier spectral methods}. The function $u(x,y,t)$ is thereby
numerically approximated via a truncated Fourier series in both $x$ 
and $y$ and \eqref{KP} is equivalent to 
\be \label{kpf}
\partial_t \widehat u + i\xi_x \widehat{u^{2}} - i \xi_x^3 \widehat u \pm \frac{i \xi_y^2}{\xi_x + i0} = 0, \quad \widehat u|_{t=0} = \widehat u_{\rm in}(\xi_x, \xi_y).
\ee 
Here, the singular multiplier $-i  \xi^{-1}_x$ is regularized in standard way as $-i (\xi_x + i 0)^{-1}$, see \cite{KSM} for more details. 
These methods show excellent convergence properties for 
smooth functions such as the ones studied in here. In addition, they have the 
advantage that only minimal numerical dissipation is introduced, 
which is important since numerical dissipation could suppress some of the 
dispersive effects we are interested in. In view of \eqref{kpf}, we see that the use of Fourier methods forces us to solve 
a large system of ordinary differential equations in 
time. Due to the third derivative in $x$ and the non-local 
term this system is in general very {\it stiff}. However, since the stiffness is 
in the linear part only, there are efficient methods to allow for time 
integration of high precision which avoids a pollution of the Fourier 
coefficients with numerical errors. In \cite{etna} and \cite{KR} it 
was shown that fourth order exponential integrators are very 
efficient for KdV-like equations and we shall use a method of Cox and 
Matthews \cite{CM} in the following. 

In the following, we shall always use $N_{x}=2^{10}$ modes in $x$-direction and $N_{y}=2^{8}$ 
modes in $y$-direction for $x\in[-\pi L_{x},\pi L_{x}]$ and 
$y\in[-\pi L_{y},\pi L_{y}]$. We choose $L_{x}=8\pi K(k)/\kappa$, 
where $K(k)$ is again the complete elliptic integral of the first kind. 
This choice ensures periodicity of the solution. For the $y$-direction we choose 
$L_{y}=2$. In all our numerical simulation, we set in (\ref{cn}) $u_{0} = x_0=0$ and $k=0.5$ so that the only remaining parameter is $\kappa \in \R$.
The elliptic functions and integrals are computed with Matlab commands to machine precision. 

To test the quality of the numerical code, we first propagate the 
cnoidal solution (\ref{cn}) with the full KP code. As mentioned 
before, the KP equation admits mass conservation, i.e., conservation of the $L^2$ norm, 
but since this feature is not implemented in the code, 
numerical errors will always lead to a time-dependent $L^2$-norm. The latter can be used as an indicator for the quality of the numerical 
results, see also the discussion in \cite{KR}. As a first test case, we first consider the KP-I equation with initial data 
$$
u_{\rm in} (x,y) = u_{\rm cn} |_{t=0} (x), \quad \mbox{and $\kappa = 2$.}
$$
We find that the $L^{\infty}$-norm of the difference 
between numerical and exact solution is of the order of 
$2.6\times 10^{-8}$ for $N_{t}=10^{4}$ time steps. In addition the numerical error in the conservation of mass  
$$\Delta:=1-\frac{L^2(t)}{L^2(0)},$$ where $L^2(t)$ is the numerically computed (discrete) $L^2$ norm at time $t\in \R$, 
is roughly $10^{-10}$. For KP-II, with the same initial data, the values are found to be almost identical. 
This shows that the quantity $\Delta$ overestimates 
numerical precision as discussed in \cite{KR} by roughly two orders 
of magnitude and can be used as a reliable indicator of numerical 
accuracy. We shall always aim for values of $\Delta$ of the order of $10^{-5}$ or 
smaller in order to ensure an accuracy much better than plotting accuracy. 
For smaller amplitudes, i.e.,  $\kappa < 2$, the accuracy is found to be even higher: For $\kappa=0.5$, we find $\Delta=10^{-12.6}$ and for 
the $L^\infty$-norm of the difference between numerical and exact 
solution $1.87\times10^{-13}$ for KP-I, and $\Delta=10^{-12.6}$ and 
$1.4\times10^{-13}$ for KP-II. Working with double precision in 
Matlab, which allows for a 
precision of $10^{-14}$, we conclude that the solution can indeed be 
propagated with machine precision.

\section{Gaussian perturbations}\label{sec:gauss}

In this section we shall consider the KP equation \eqref{KP} subject to initial data of the form
\begin{equation}\label{pini}
u_{\rm in} (x,y) = u_{\rm cn} |_{t=0}(x) + u_{\rm p}(x,y),
\end{equation}
where $u_{\rm cn}$ is the cnoidal wave for $t=0$ (and $u_0 = x_0 = 
0$) and $u_{\rm p}$ is a (small) perturbation which is assumed to be {\it localized} with respect to both $x$ and $y$. 

Next, we recall the well known fact (see for instance \cite{MST} for a rigorous proof) that KP equations satisfy the following {\it constraint}, in the sense of Riemann Integrals:
\begin{equation}
    \int_{X} u (t,x,y) \, dx=0, \quad y \in X, t\not =0,
    \label{constraint}.
\end{equation}
where $X= \R$, or $X=\mathbb T$, respectively.
Equation \eqref{constraint} holds, even if the initial condition $u_{\rm in}$ does 
not satisfy the constraint, see \cite{KSM} for a numerical example. 
In such a situation the solution ceases to be continuous at $t=0$, however. Since this leads to numerical problems, we will 
only consider initial conditions which satisfy the above constraint. This is 
obviously the case for the cnoidal KdV solution, but we shall also 
impose it on the perturbations $u_{\rm p}$. A possible way to do so is to 
consider $u_{\rm in} = \partial_x f(x,y)$, where $f$ is a periodic or Schwarzian function, which 
can be considered as essentially periodic if the period is chosen 
large enough such that
the solution at the limits of the computational domain decreases 
below machine precision. The latter is important to avoid the appearance of {\it Gibbs phenomena} 
due to discontinuities.

\begin{figure}[htb]
    \begin{center}
      \includegraphics[width=\textwidth]{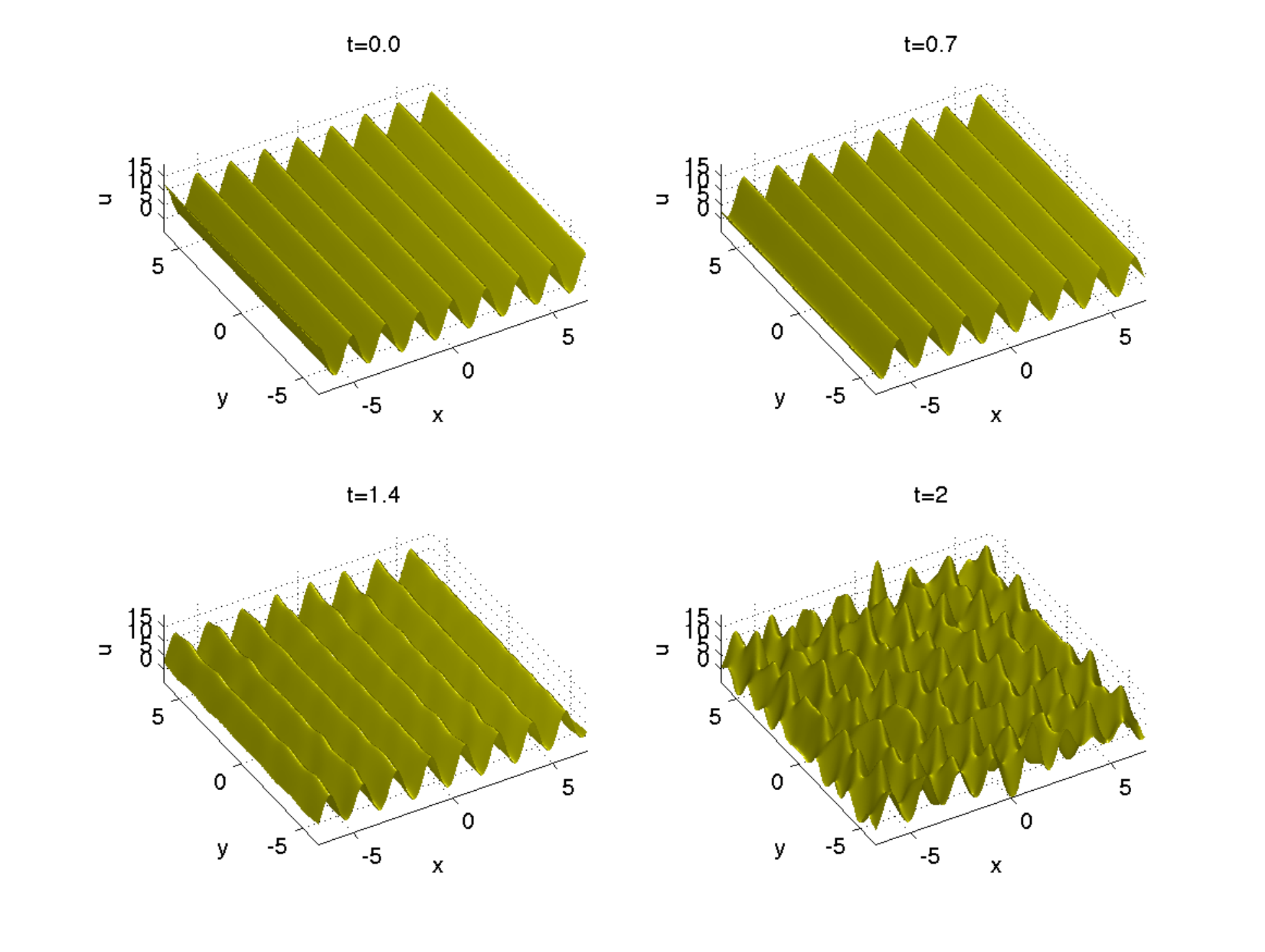}
    \end{center}
\caption{Solution to the KP-I equation (\ref{KP}) for several values of $t$ obtained for initial data given by a cnoidal wave (\ref{cn}) with parameters
$\kappa=2$, $k=0.5$, and $x_0=u_{0}=0$, and a Gaussian perturbation $u_{\rm p}(x,y)$ of the form
(\ref{pert}).}
\label{kpI2gauss}
\end{figure}

\begin{figure}[htb!]
 \centering
 \includegraphics[width=0.6\textwidth]{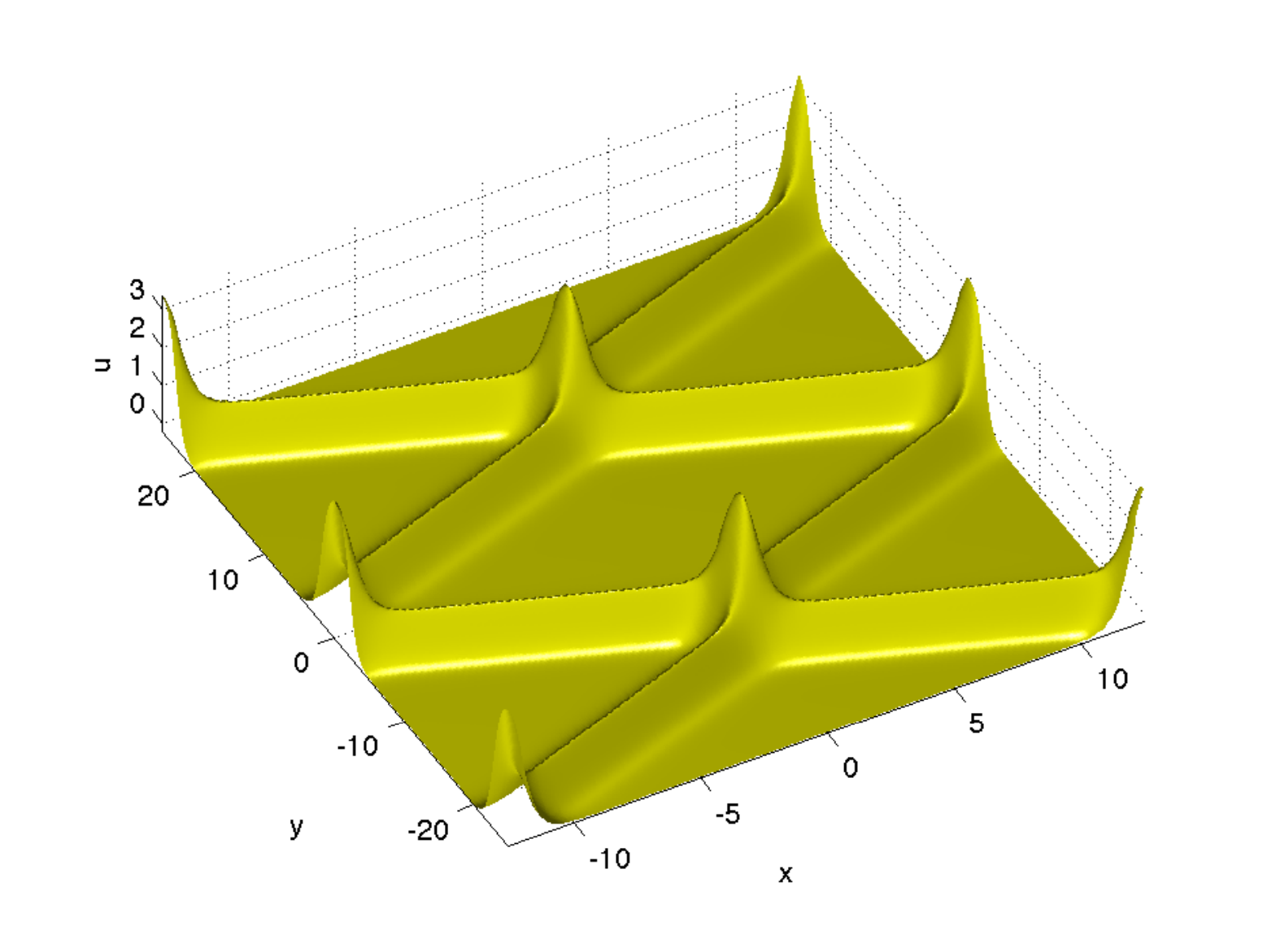}
 \caption{Doubly periodic hyperelliptic KP-II solution of genus 2 for 
 a fixed time. The profile travels with constant speed in 
 $x$-direction.}
 \label{hyper}
\end{figure}

To be more precise, we shall study localized perturbations of the following form 
\begin{equation}
    u_{\rm p}(x,y)=xe^{-(x^{2}+y^{2})}
    \label{pert},
\end{equation}
i.e., an $x$-derivative of a Gaussian. With the same parameters as 
above and $\kappa=2$, we obtain for the time evolution of the initial 
data \eqref{pini} that at $t=2$ the relative $L^2$ norm is $\Delta\sim 10^{-10}$. It can be seen in Fig.~\ref{kpI2gauss} that the 
solution stays close to the cnoidal form for some time, but will 
finally disintegrate into a periodic array of humps. 
It is conjectured that for $(x,y)\in \mathbb{R}^{2}$, solutions to KP-I 
asymptotically decompose into several \emph{lump solitons}, i.e., solitons which are 
localized in both spatial directions with an algebraic decay towards infinity. Numerically, this was already observed in \cite{KS}. 
However, there is no periodic analogue of the lump soliton known, and 
thus it is not clear whether there are 
doubly periodic solutions to KP-I which appear as asymptotic 
solutions of a perturbed cnoidal solution. There are, however, doubly periodic solutions 
to the KP equations, which can be given in terms of multi-dimensional 
theta functions on hyperelliptic Riemann surfaces, see, e.g., 
\cite{dubrovin, lmp}. The simplest such example is a travelling wave 
on a genus 2 surface that, which for KP-II is depicted in Fig.~\ref{hyper}.

\begin{figure}[htb]
    \begin{center}
      \includegraphics[width=\textwidth]{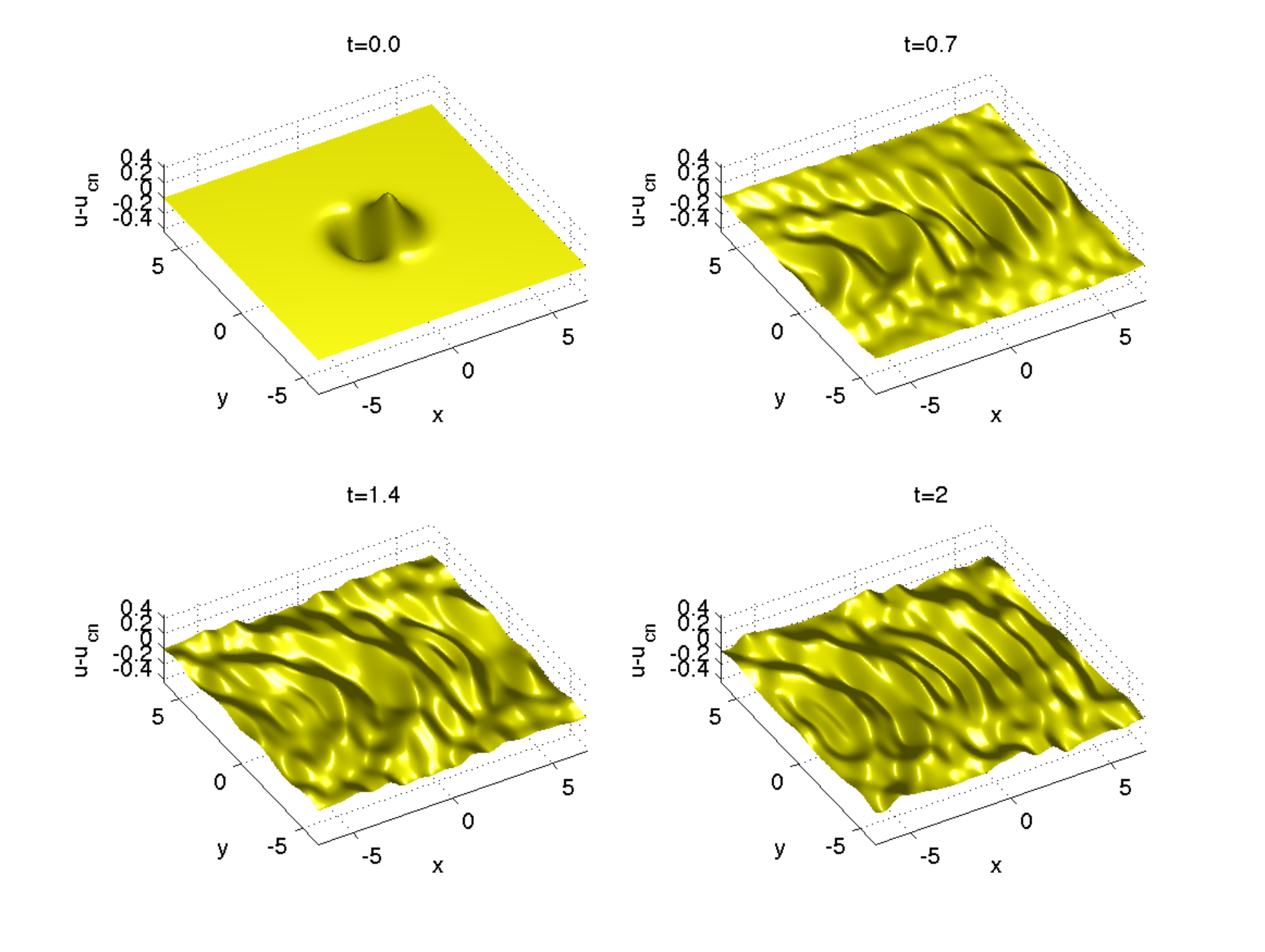}
    \end{center}
\caption{Difference of the solution to the KP-II equation (\ref{KP}) for the initial 
data \eqref{pini}  with $u_{\rm p}$ given by(\ref{pert}), respectively. 
The cnoidal initial data are chosen with  
$\kappa=2>\kappa_c$, $k=0.5$ and $x_0=u_{0}=0$.}
\label{kpII2gauss}
\end{figure}

The situation is considerably different if the same initial data \eqref{pert} are 
used for the KP-II equation (for which we essentially obtain the same value of $\Delta$). Note that even though the equation is defocusing, 
the perturbation cannot be radiated away to infinity as would be the case for the perturbations 
of the KdV soliton on all of $\mathbb{R}^{2}$. In our periodic setting, they 
can just smooth out over the whole period or distribute periodically. 
In Fig.~\ref{kpII2gauss} we show the difference between the KP-II 
solution corresponding to initial data \eqref{pini} with $u_{\rm p} =0$ (i.e. a purely cnoidal solution) and $u_{\rm p}$ given by \eqref{pert}, respectively. It can be 
seen that, as time increases, the initial hump in fact smoothes itself out over the whole domain. 

\begin{figure}[htb]
    \begin{center}
      \includegraphics[width=\textwidth]{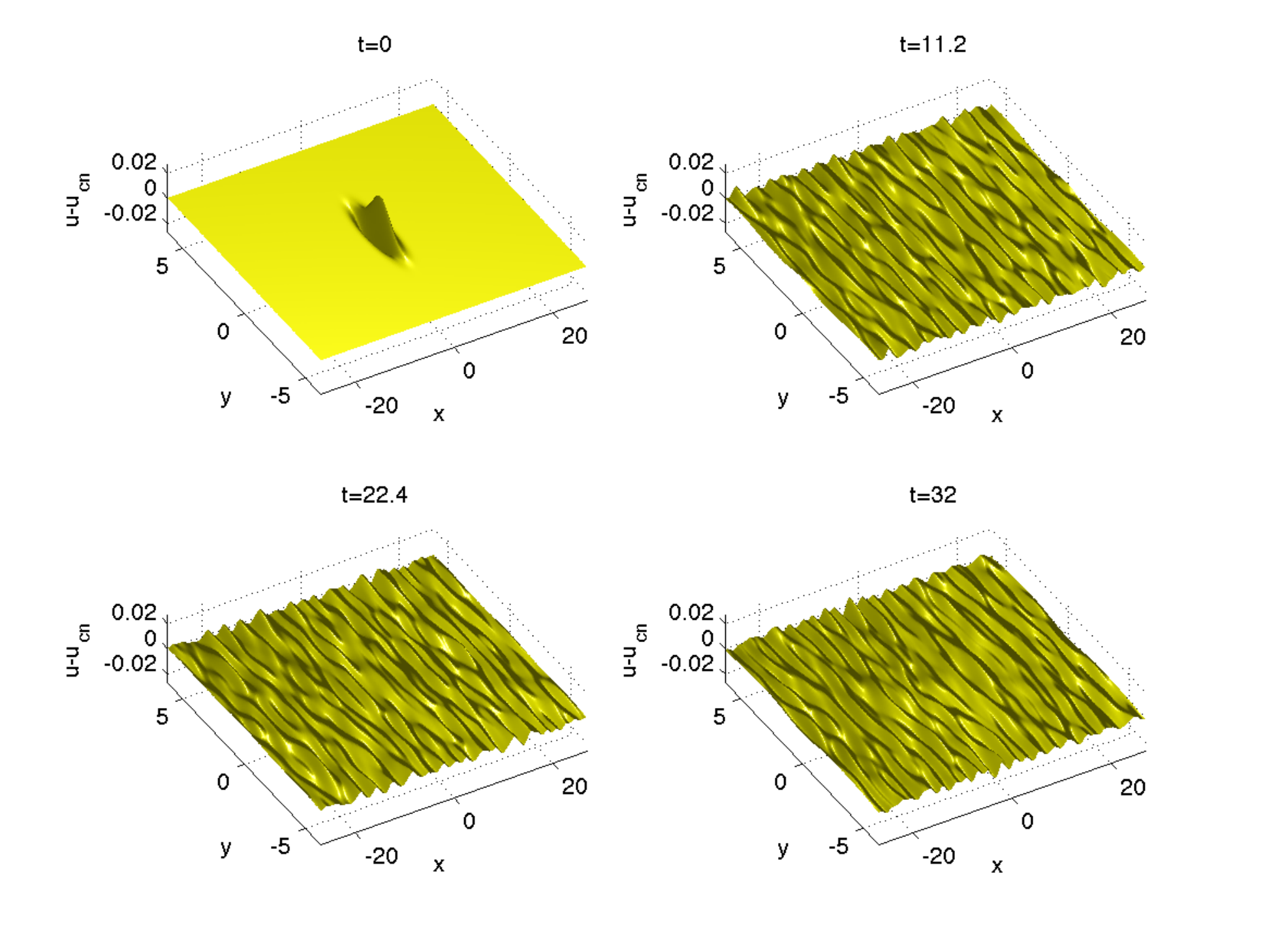}
    \end{center}
\caption{Difference of the solution to the KP-II equation (\ref{KP}) for initial 
data $u_{\rm in}=u_{\rm cn}+\frac{1}{16}u_{\rm p}$ with  $u_{\rm p}$ given by
(\ref{pert}), respectively. The cnoidal initial data are chosen with  
$\kappa=0.5<\kappa_c$, $k=0.5$ and $x_0=u_{0}=0$.}
\label{kpII05gaussdiff}
\end{figure}

So far, all simulations were done with $\kappa = 2$ which is above the critical value $\kappa_c$ established in \cite{RT2} for the 
transverse stability of solitary wave (obtained for $k\to 1$ in the cnoidal solution). Our simulations confirm that for $\kappa>\kappa_{c}$ 
the cnoidal solution is unstable for KP-I, but stable for KP-II. To check whether there is a similar 
threshold below which cnoidal waves become stable for KP-I, we shall also consider the value  
$\kappa=0.5$ which is well below $\kappa_{c}$. To compare similar situations, 
we multiply the final time by a factor of 16 (since the the velocity $V$ in 
(\ref{cn}) scales as $\kappa^{2}$) and in addition consider a perturbation 
of the form (\ref{pert}) but rescaled by a factor $u_{p}/16$. The latter ensures that the perturbation is of the same relative 
amplitude as in the case $\kappa=2$ considered above. The code produces for both KP-I and KP-II an error $\Delta<10^{-12}$. 
In Fig.~\ref{kpII05gaussdiff} it can be seen that the KP-II 
situation is  essentially unchanged in this case, i.e., the cnoidal 
solution is stable for KP-II also for $\kappa < \kappa_c$. 

\begin{figure}[htb]
    \begin{center}
      \includegraphics[width=\textwidth]{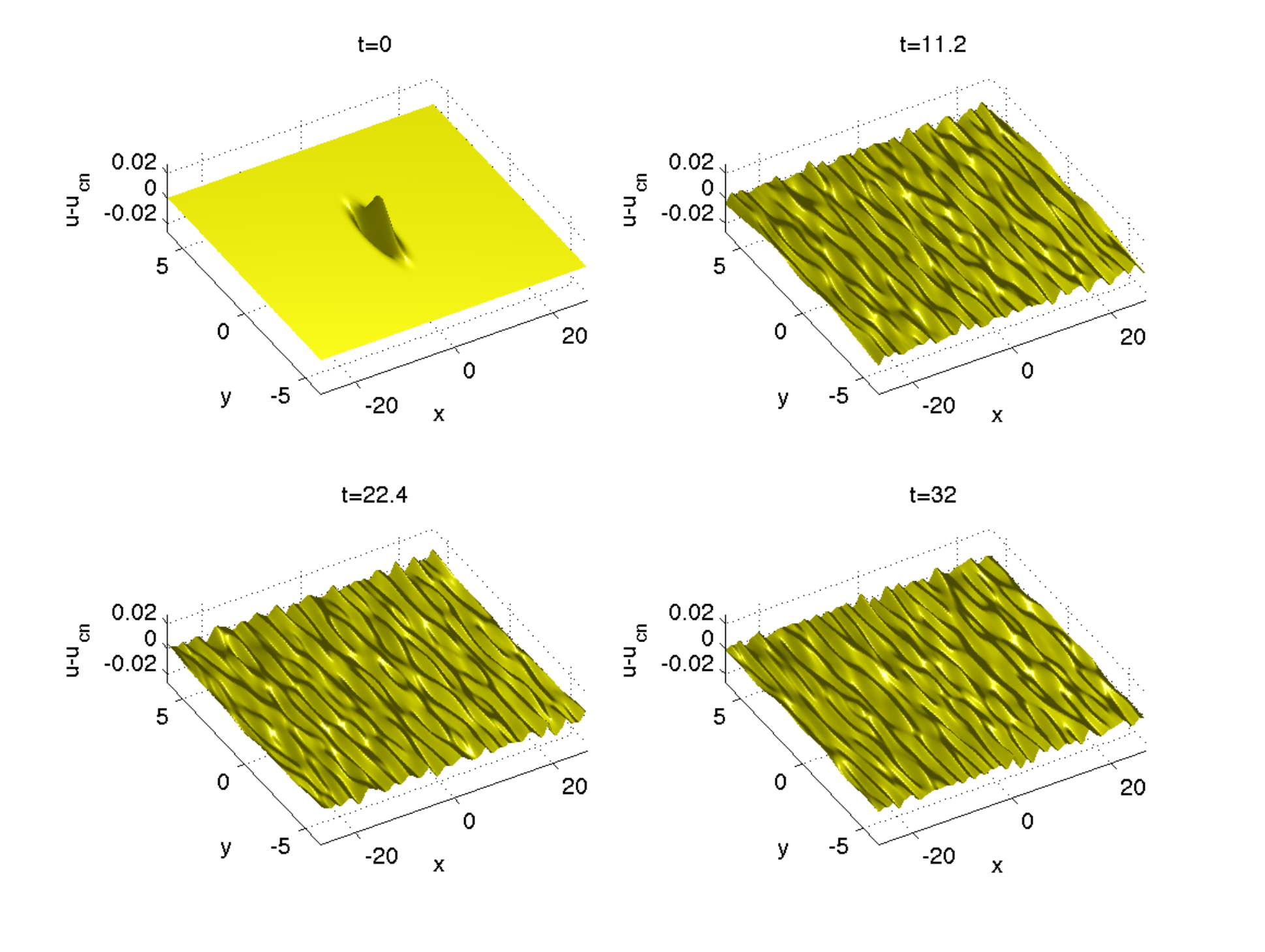}
    \end{center}
\caption{Difference of the solution to the KP-I equation (\ref{KP}) for initial 
data $u_{\rm in}=u_{\rm cn}+\frac{1}{16}u_{\rm p}$ with  $u_{\rm p}$ given by
(\ref{pert}), respectively. The cnoidal initial data are chosen with  
$\kappa=0.5<\kappa_c$, $k=0.5$ and $x_0=u_{0}=0$.}
\label{kpI05gauss}
\end{figure}

The corresponding situation for KP-I is found to be considerably different when compared to the situation shown in Fig.~\ref{kpI2gauss}. 
Indeed, it can be seen in Fig.~\ref{kpI05gauss} that for $\kappa=0.5<\kappa_c$ the solution is almost identical to the KP-II case, and thus also stable. 
It consequently appears that, just as for solitary waves, there is also a critical value for $\kappa$ 
separating the stable from the unstable regime in KP-I.

\section{Periodic deformations}\label{sec:per}
In this section we shall consider initial perturbations of \eqref{cn}, which are {\it no longer spatially localized}. Rather, we shall consider a $y$-periodic deformation of the 
of the form
\begin{equation}
    u_{\rm in}(x,y) =12\kappa^{2}k^{2}\mbox{cn}^{2}\left(k\left(x+\delta \cos 
    \frac{4y}{L_{y}}\right);k\right),\quad \delta \ge0.
    \label{def},
\end{equation}
For $\delta = 0$ this yields the usual cnoidal wave (with $x_0=u_0=0$). In the following, though, we shall choose $\delta=0.4$ yielding a 
strongly deformed cnoidal type solution as initial data. Note that these kind of $y$-oscillatory perturbations are different from the $x$-periodic 
perturbations considered in, e.g., \cite{HSS, Ha}.

\begin{figure}[htb]
    \begin{center}
      \includegraphics[width=\textwidth]{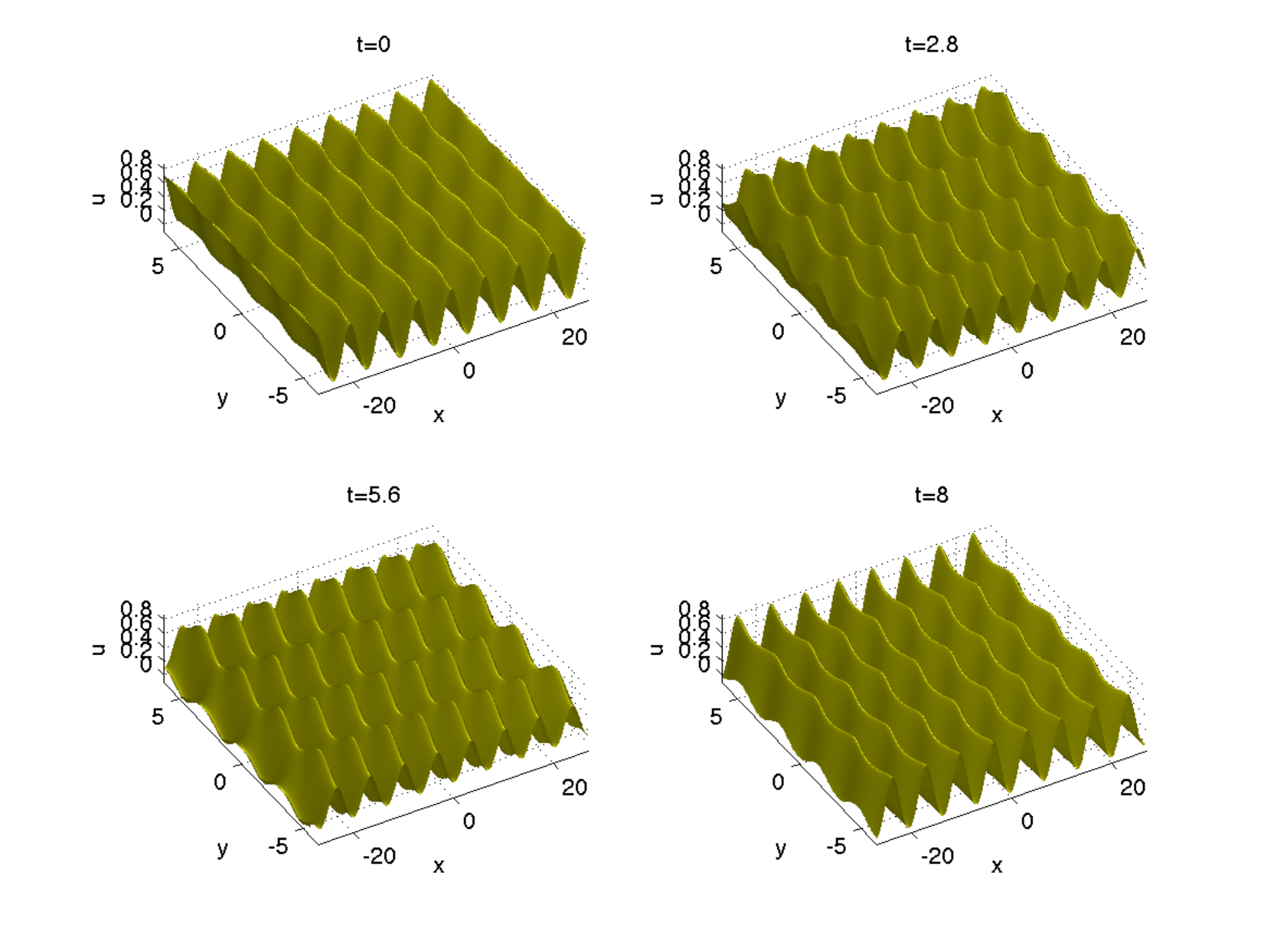}
    \end{center}
\caption{Solution to the KP-II equation for several values of $t$, obtained for initial 
data given by (\ref{def}), i.e. periodically deformed cnoidal initial data with 
$\kappa=0.5$, $k=0.5$ and $x_0=u_{0}=0$.}
\label{kpII05cos}
\end{figure}

\begin{figure}[htb]
    \begin{center}
	\includegraphics[width=0.45\textwidth]{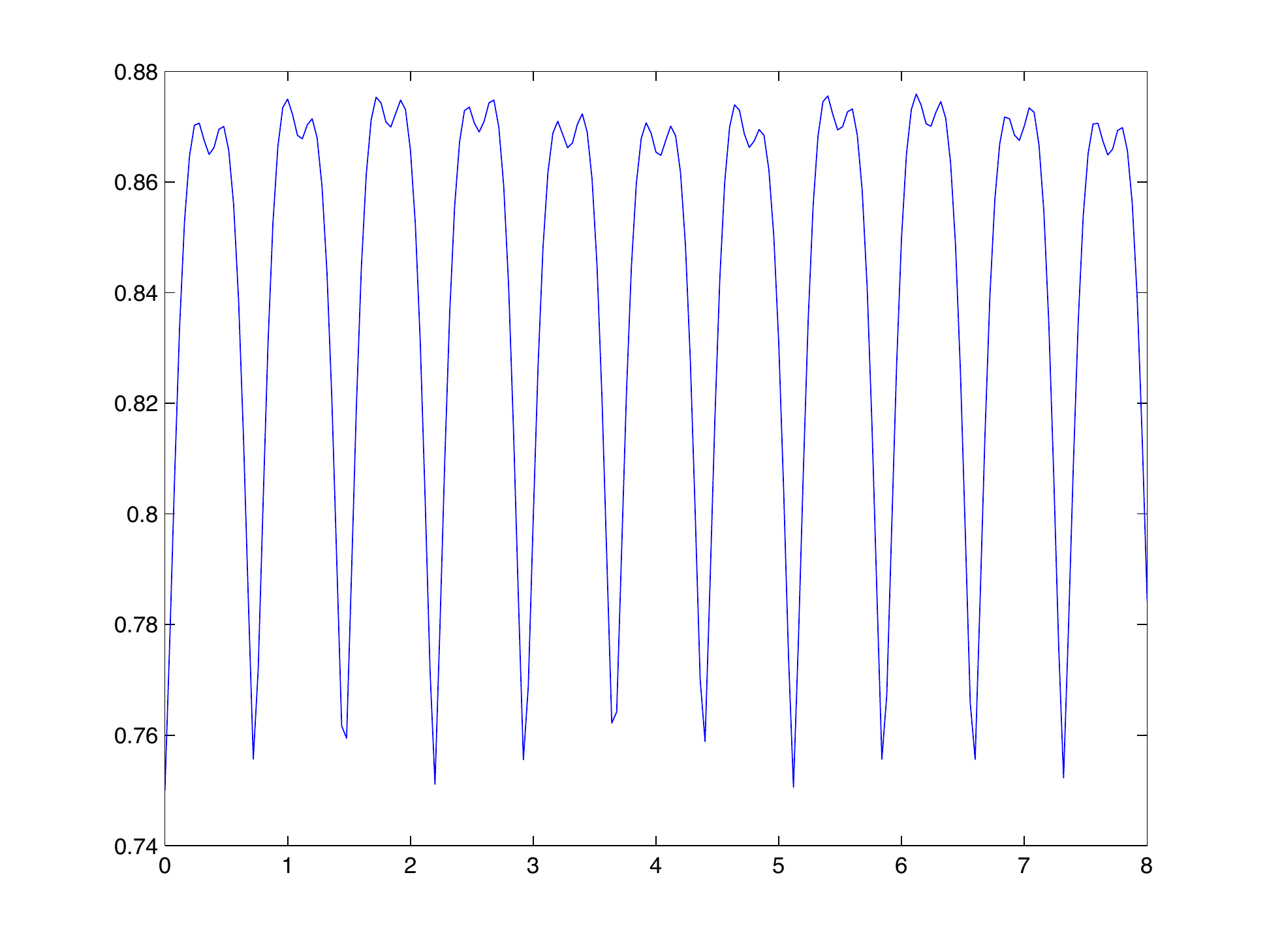}
	\includegraphics[width=0.45\textwidth]{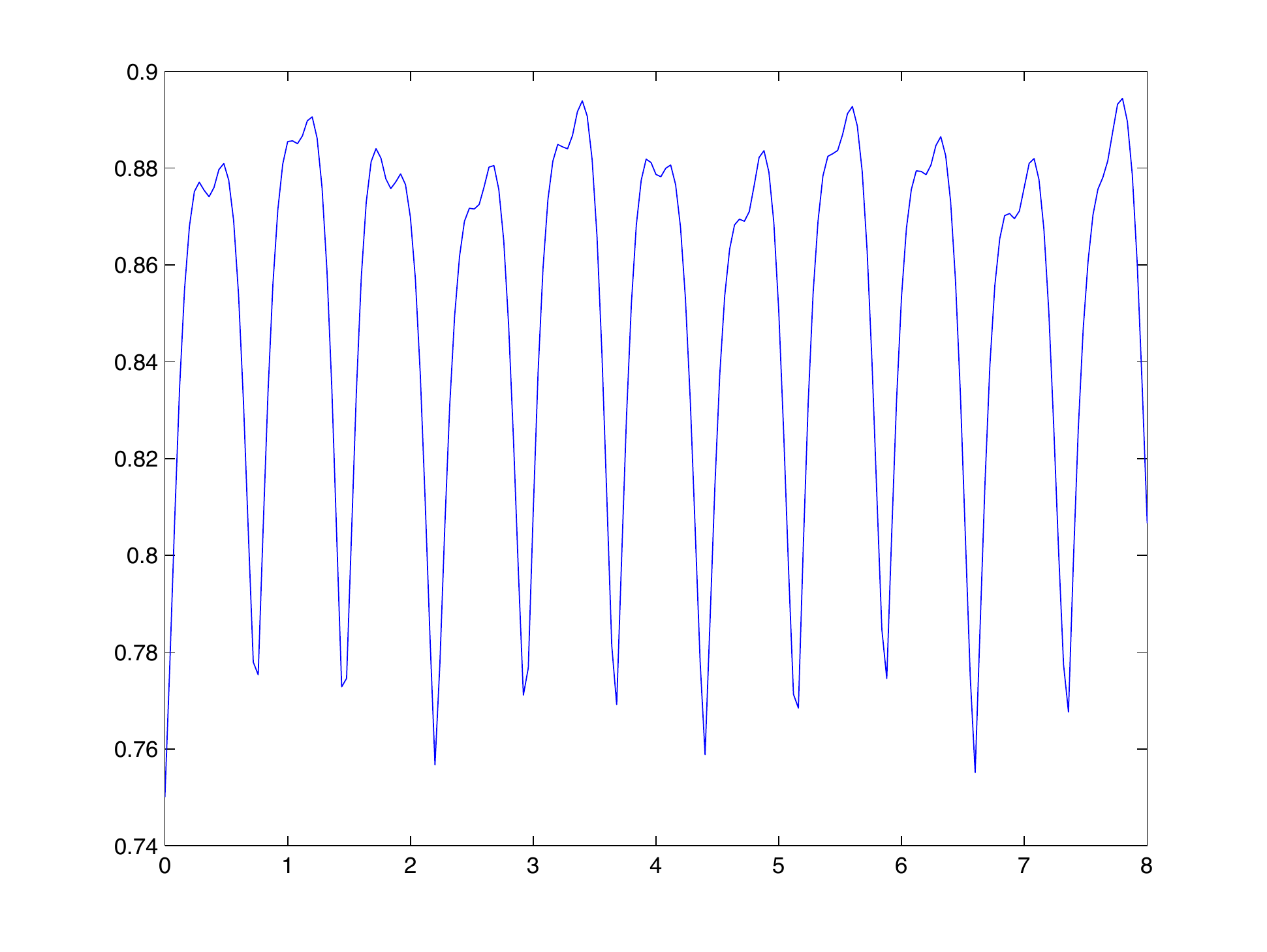}
    \end{center}
\caption{$L^\infty$ norm of the KP-II solution shown in 
Fig.~\ref{kpII05cos} as a function on time on the left and for KP-I 
of Fig.~\ref{kpI05cos} on the right.}
\label{kpII05cosmax}
\end{figure}

\begin{figure}[htb]
    \begin{center}
      \includegraphics[width=\textwidth]{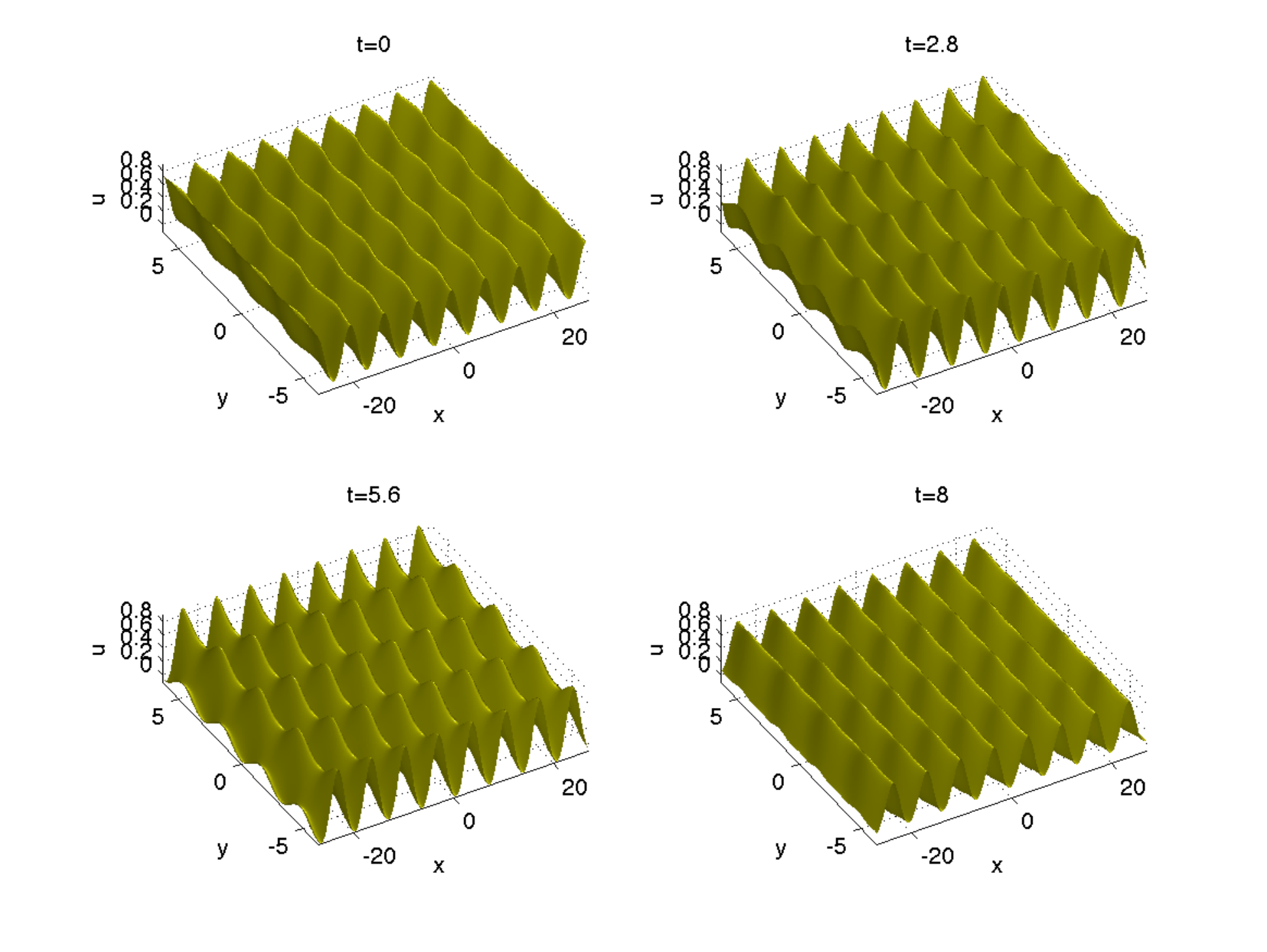}
    \end{center}
\caption{Solution to the KP-I equation for several values of $t$, obtained for initial 
data given by (\ref{def}), i.e. periodically deformed cnoidal initial data with 
$\kappa=0.5$, $k=0.5$ and $x_0=u_{0}=0$.}
\label{kpI05cos}
\end{figure}

For small $\kappa =0.5$ these initial data lead a solution of KP-II equation as shown in Fig.~\ref{kpII05cos}.
In fact the solution oscillates around the cnoidal solution which 
thus appears stable in this sense. Indeed, these oscillations can be seen best 
by tracing the $L^\infty$ norm of the solution as a function 
of time, see Fig.~\ref{kpII05cosmax}. One can expect these 
oscillations to go on indefinitely.
For the KP-I, the same initial data lead to an essentially identical behavior as can be seen in Fig.~\ref{kpII05cosmax} and 
Fig.~\ref{kpI05cos}. Thus for small $\kappa$, KP-I and KP-II show 
very similar behavior for perturbed cnoidal initial data.

\begin{figure}[htb]
    \begin{center}
      \includegraphics[width=\textwidth]{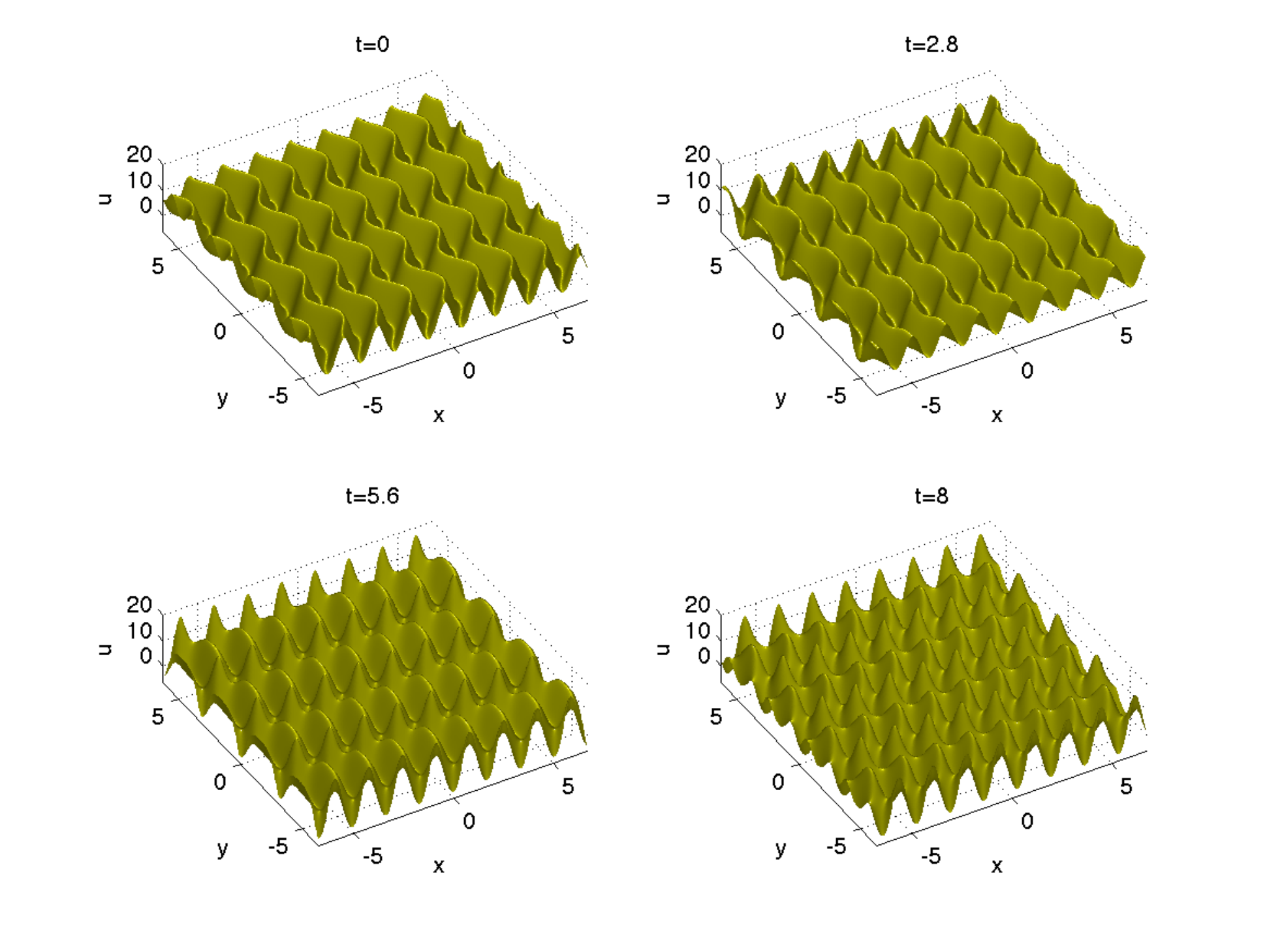}
    \end{center}
\caption{Solution to the KP-I equation (\ref{KP}) as a function of time, obtained from initial data given by (\ref{def}) i.e., periodically deformed cnoidal initial data with 
$\kappa=2$, $k=0.5$ and $x_0=u_{0}=0$.}
\label{kpI2cos}
\end{figure}

Analogous to the case of localized perturbations, the situation becomes different for larger $\kappa$: For $\kappa=2$, the KP-I solution corresponding to initial data (\ref{def}) 
is given in Fig.~\ref{kpI2cos}. The solution appears to be periodic in time in 
the sense that there are times where it is close to the initial 
configuration. But there is much more structure in this type of {\it doubly 
periodic solution} than in the genus 2 example shown in Fig.~\ref{hyper}. In 
particular there seems to be a `breathing' structure similar to the so-called {\it breathers} in the {\it nonlinear Schr\"odinger equation} (NLS), see, e.g., \cite{DT} 
(for the relation between NLS and KP solutions see \cite{CR, kalla}).

\begin{figure}[htb]
    \begin{center}
      \includegraphics[width=\textwidth]{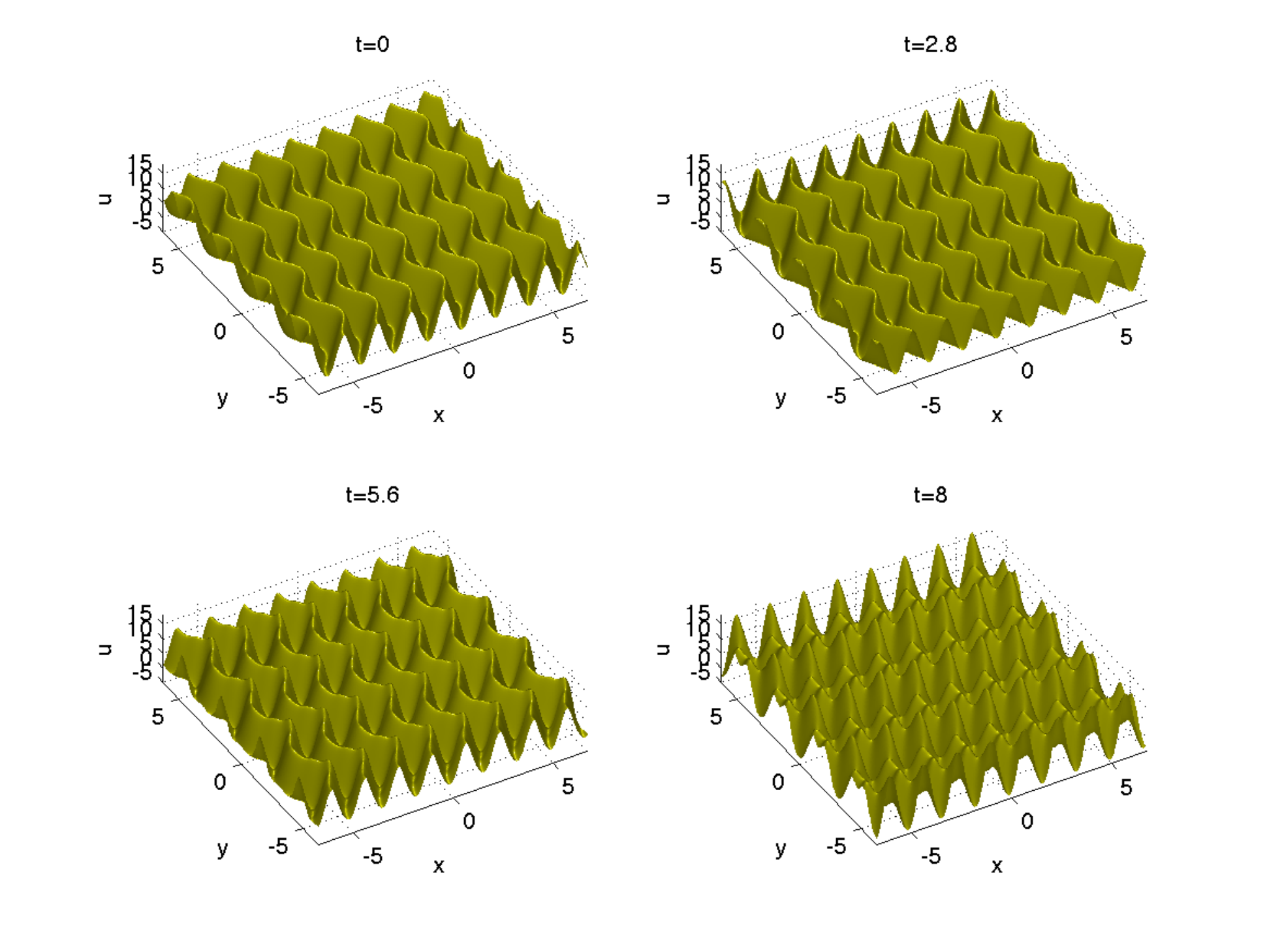}
    \end{center}
\caption{Solution to the KP-II equation (\ref{KP}) as a function of time, obtained 
for initial data given by (\ref{def}), i.e. periodically deformed cnoidal initial data with 
$\kappa=2$, $k=0.5$ and $x_0=u_{0}=0$.}
\label{kpII2cos}
\end{figure}

For KP-II the situation for $\kappa=2$ is similar to the one of KP-I as can be seen in 
Fig.~\ref{kpII2cos}. Again, the initial configuration reappears from time to 
time, but there are very rich oscillations patterns. 
It thus seems to be that for such kind of de-localized ($y$-periodic) perturbations, both KP-I and KP-II exhibit a similar behavior and both 
are not stable in the usual sense. 

\section{Conclusion}
We presented numerical simulations concerning the transverse 
stability and instability of cnoidal waves for KP equations. By means of a spectral method we found that:
\begin{itemize}
\item For spatially localized perturbations KP-II is stable for general $\kappa \in \R$, whereas KP-I is only stable for $\kappa$ below a certain threshold (and unstable above). 
This is in agreement with recent rigorous results on the transverse stability of solitary waves (obtained from the cnoidal solution in the limit $k \to 1$).
\item For $y$-periodic perturbations and small $\kappa$, both, KP-I and KP-II are stable. For large $\kappa$ both KP-I and KP-II appear to be unstable in the usual sense, as both 
exhibit a doubly periodic type solution with rich oscillatory behavior. It remains a challenging open problem to give a more precise description of 
the long time behavior of cnoidal waves under such type of deformations.
\end{itemize}

\bibliographystyle{amsalpha}

\end{document}